\documentclass[12pt]{article}
\usepackage{amsmath}
\usepackage{amssymb}
\usepackage{amsmath,amssymb,amsthm,amsfonts}
\newtheorem{theorem}{Theorem}
\newtheorem{lemma}{Lemma}
\newtheorem{corollary}{Corollary}

{\theoremstyle{definition}

\newtheorem{remark}{Remark}
\newtheorem{example}{Example}
}

\newcommand{\il}[2]{\int\limits_{#1}^{#2}}

\newcommand{\ph}{\phantom{a}}
\newcommand{\phh}{\phantom{aaa}}

\begin{document}

MSC 34D20

\vskip 20pt

\centerline{\bf A generalization of the Gronwall-Bellman lemma}

\vskip 10 pt

\centerline{\bf G. A. Grigorian}

\vskip 10 pt

\centerline{0019 Armenia c. Yerevan, str. M. Bagramian 24/5}
\centerline{Institute of Mathematics of NAS of Armenia}
\centerline{E - mail: mathphys2@instmath.sci.am, \ph phone: 098 62 03 05, \ph 010 35 48 61}

\vskip 20 pt

\noindent
Abstract. The Riccati equation method is used to obtain a generalization of the Gronvall-bellman lemma te obtained result is used to generalize a result of Lyapunov.

\vskip 20 pt

\noindent
Key words:  the Gronwall-Bellman lemma, the Riccati equation method, the Lyapunov estimate. 

\vskip 20 pt

{\bf  1. Introduction}. The Gronwal-Bellman lemma ([3, p. 35],  \ph [7. p. 108])  states

\begin{lemma} (Gronwall-Bellman). Let $u(t)$ and $v(t)$ be nonnegative continuous functions on $[t_0,+\infty)$. If
$$
u(t) \le c + \il{t_0}{t}v(\tau)u(\tau) d \tau, \phh t \ge t_0,
$$
where $c$ is a positive constant, then
$$
u(t) \le c\exp\bigg\{\il{t_0}{t}v(\zeta) d\zeta\biggr\}, \phh t \ge t_0.
$$
\end{lemma}
\phantom{aaaaaaaaaaaaaaaaaaaaaaaaaaaaaaaaaaaaaaaaaaaaaaaaaaaaaaaaaaaaaaaaaaaaa} $\blacksquare$

This lemma serves as a important tool for qualitative study of ordinary differential equations and many generalizations in different directions of it have been obtained (see e.g. [2, 4--6, 8, 10--22] and cited works therein).

In this short  note we use a comparison criteria for scalar  Riccati equations to prove the following generalization of Lemma 1.

\begin{lemma}
Let $u(t), \ph v(t)$ and $f(t)$ be real-valued continuous functions on $[t_0,+\infty)$ such that
$$
u(t) \ge 0, \ph v(t)\ge 0,  \ph   \il{t_0}{t}f(\zeta) d\zeta \ge 0,
$$
$$
u(t) \le c + \il{t_0}{t}v(\zeta) u(\zeta) d \zeta + \il{t_0}{t}f(\zeta) d\zeta, \ph t \ge t_0, \eqno (1.1)
$$
where $c$ is  a positive constant, then
$$
u(t) \le c\exp\biggl\{\il{t_0}{t}v(\zeta) d\zeta\biggr\} + \il{t_0}{t}\exp\biggl\{\il{\zeta}{t}v(s) d s\biggr\}f(\zeta) d \zeta, \phh t \ge t_0. \eqno (1.2)
$$
\end{lemma}
\phantom{aaaaaaaaaaaaaaaaaaaaaaaaaaaaaaaaaaaaaaaaaaaaaaaaaaaaaaaaaaaaaaaaaaaaaa} $\blacksquare$

\begin{corollary}
Let a positive continuous on $(a,b)$ function $u(t)$ for every $t, \ph  \tau \in (a,b)$ satisfy the inequality
$$
u(t) \le u(\tau) + \il{\tau}{t} v(\zeta) u(\zeta) |d\zeta| + \il{\tau}{t} f(\zeta)|d \zeta|, \eqno (1.3)
$$
where $v(t), \ph f(t) \in \mathbb{C}(a,b), \ph v(t) \ge 0$ and $f(t) \ge 0$ for $a < t < b.$ Then for  $a < t_0  \le t < b$ the following inequalities are valid
$$
u(t_0)\exp\biggl\{-\il{t_0}{t}v(\zeta) d \zeta\biggr\} - \il{t_0}{t}\exp\biggl\{-\il{\zeta}{t}v(s) d s\biggr\} f(\zeta) d \zeta \le u(t) \le  \phantom{aaaaaaaaaaaaaaaaaaaaaaaaaaaaaaaa}
$$
$$
\phantom{aaaaaaaaaaaaa}\le u(t_0)\exp\biggl\{\il{t_0}{t}v(\zeta) d \zeta\biggr\} + \il{t_0}{t}\exp\biggl\{\il{\zeta}{t}v(s) d s\biggr\} f(\zeta) d \zeta. \eqno (1.4)
$$
\end{corollary}
\phantom{aaaaaaaaaaaaaaaaaaaaaaaaaaaaaaaaaaaaaaaaaaaaaaaaaaaaaaaaaaaaaaaaaaaaaa} $\blacksquare$

\begin{example}
Consider the linear system of ordinary differential equations
$$
Y' = A(t) Y + g(t), \ph t \ge t_0,
$$
where $A(t)$ is a real-valued continuous matrix function of dimension $n\times n$ on $[t_0,+\infty), \ph g(t)$ is a real-valued continuous vector function of dimension $n$ on $[t_0,+\infty)$. Let $Y(t)$ be a solution of this equation. Then
$$
Y(t) = Y(\tau) +\il{\tau}{t}A(\zeta)Y(\zeta) d\zeta + \il{\tau}{t}g(\zeta)d\zeta, \phh t, \tau \ge t_0.
$$
From here we obtain
$$
||Y(t)|| \le ||Y(\tau)|| + \il{\tau}{t}||A(\zeta)|| ||Y(\zeta)|| d\zeta + \il{\tau}{t}||g(\zeta)||d\zeta, \phh t, \tau \ge t_0,
$$
where $||Y(t)||$ and $||g(t)||$ denote  euclidian norms of $Y(t)$ and $g(t)$ respectively   in $\mathbb{R}^n$, and $||A(t)||$ denotes the norm of the operator $A(t) : \mathbb{R}^n \to \mathbb{R}^n$ for fixed $t\ge t_0$. Then by Corollary 1 we have the following estimates
$$
||Y(t_0)||\le \exp\biggl\{-\il{t_0}{t}||A(\zeta)||d\zeta\biggr\} - \il{t_0}{t}\exp\biggl\{-\il{\zeta}{t}||A(s)||ds\biggr\}||g(\zeta)|| d\zeta \le \phantom{aaaaaaaaaaaaaaaa}
$$
$$
\phantom{aaaaaaaaaaaaaaaa} \le ||Y(t)|| \le \exp\biggl\{\il{t_0}{t}||A(\zeta)||d\zeta\biggr\} + \il{t_0}{t}\exp\biggl\{-\il{\zeta}{t}||A(s)||ds\biggr\}||g(\zeta)|| d\zeta
$$
Note that these estimates generalize the Lyapunov's result (see [1, p. 132])
\end{example}

{\bf 2. Auxiliary propositions}. Let $f_k(t), \ph g_k(t), \ph h_k(t), \ph k=1,2,$ be real-valued continuous functions on $[t_0,+\infty)$. Consider the Riccati equations
$$
y' + f_k(t) y^2 + g_k(t) y + h_k(t) = 0, \phh t\ge t_0, \eqno (2.4_k)
$$
$k=1,2$ and the differential inequalities
$$
\eta' + f_k(t) \eta^2 + g_k(t) \eta + h_k(t) \ge 0, \phh t\ge t_0, \eqno (2.5_k)
$$
$k=1,2$.

\begin{remark}
\label{Remark 2.1.} Every solution  of Eq. $(2.4_2)$ on $[t_0,t_1)$ is also a solution of the inequality $(2.5_2)$ on $[t_0,t_1)$.
\end{remark}

\begin{remark}\label{Remark 2.2.} If $f_1(t) \ge 0, \ph t\in [t_0,t_1)$, then every solution of the linear equation
$$
\zeta' + g_1(t)\zeta + h_1(t) = 0, \phh t\in [t_0,t_1)
$$
is also a solution of the inequality $(2.5_1)$ on $[t_0,t_1)$.
\end{remark}

\begin{theorem}\label{Theorem 2.1} [9, Theorem 3.1].  Let $y_2(t)$ be a solution of Eq. $(2.4_2)$ on $[t_0,\tau_0) \linebreak (t_0 < \tau_0 \le +\infty)$ and let $\eta_1(t)$ and $\eta_2(t)$ be solutions of the inequalities $(2.5_1)$ and $(2.5_2)$  respectively on $[t_0,\tau_0)$ such that $y_2(t_0) \le \eta_k(t_0) \ph k=1,2.$ In addition let the following conditions be satisfied: $f_1(t) \ge 0, \ph \gamma - y_2(t_0) + \il{t_0}{t}\exp\biggl\{\il{t_0}{\tau}[f_1(s)(\eta_1(s) + \eta_2(s)) + g_1(s)]ds\biggr\}\biggl[(f_2(\tau) - f_1(\tau))^2 y_2^2(\tau) + (g_2(\tau) - g_1(\tau)) y_2(\tau) + h_2(\tau) - h_1(\tau)\biggr] d \tau \ge 0, \ph t\in [t_0,\tau_0)$ for some $\gamma \in [y_2(t_0), \eta_1(t_0)]$. Then Eq. $(2.4_1)$ has a solution $y_1(t)$ on $[t_0,\tau_0)$ with $y_1(t_0) \ge \gamma$ and $y_1(t) \ge y_2(t), \ph t\in [t_0,\tau_0)$.
\end{theorem}
\phantom{aaaaaaaaaaaaaaaaaaaaaaaaaaaaaaaaaaaaaaaaaaaaaaaaaaaaaaaaaaaaaaaaaaaaaaaa}$\blacksquare$

\vskip 10pt

{\bf 3. Proof of Lemma 2.} It follows from the conditions of the lemma that
$$
c + \il{t_0}{t}v(\zeta)u(\zeta) d\zeta > 0 \phh t \ge t_0.
$$
Then the function
$$
\mathcal{A}(t)\equiv \frac{v(t)\Bigl[c + \il{t_0}{t}v(\zeta)u(\zeta) d\zeta - u(t) + F(t)\Bigr]}{\Bigl[c + \il{t_0}{t}v(\zeta)u(\zeta) d\zeta \Bigr]^2}, \phh t \ge t_0
$$
where $F(t) \equiv \il{t_0}{t}f(\zeta)d\zeta, \phh t \ge t_0$ is defined on $[t_0,+\infty)$. It follows from the conditions of the lemma that
$$
\mathcal{A}(t) \ge 0, \phh t \ge t_0. \eqno (3.1)
$$
Consider the Riccati equation
$$
y' + \mathcal{A}(t) y^2 = v(t) y + v(t)F(t), \phh t \ge t_0. \eqno (3.2)
$$
and the linear one
$$
x' = U(t)x + v(t) F(t), \phh t \ge t_0. \eqno (3.3)
$$
Obviously, $y(t)\equiv c + \il{t_0}{t}v(\zeta)u(\zeta) d\zeta, \ph t \ge t_0$ is a solution of Eq. (3.2) on $[t_0,+\infty)$ and according to the Cauchy formula the function
$$
x(t)\equiv \exp\biggl\{\il{t_0}{t}v(\zeta)d\zeta\biggr\}\biggl[c + \il{t_0}{t}\exp\biggl\{-\il{t_0}{t}v(s)ds\biggr\}v(\zeta)F(\zeta)d\zeta\biggr], \phh t \ge t_0
$$
is a solution of Eq. (3.3) on $[t_0,+\infty)$. We can interpret Eq. (3.3) as a Riccati equation with $\equiv 0$ coefficient at $x^2$. Then applying Theorem \ref{Theorem 2.1} to the pair of equations (3.2) and (3.3) and taking into account (3.1) we get $y(t) \le x(t), \ph t \ge t_0$, i. e.
$$
c + \il{t_0}{t}v(\zeta)u(\zeta) d\zeta \le \exp\biggl\{\il{t_0}{t}v(\zeta)d\zeta\biggr\}\biggl[c + \il{t_0}{t}\exp\biggl\{-\il{t_0}{t}v(s)ds\biggr\}v(\zeta)F(\zeta)d\zeta\biggr], \eqno (3.4)
$$
$t \ge t_0.$. Using the integration by parts rule we will have
$$
\il{t_0}{t}\exp\biggl\{-\il{t_0}{\zeta}v(s) d s\biggr\}v(\zeta) F(\zeta) d\zeta = - \il{t_0}{t}\biggl[\exp\biggl\{-\il{t_0}{\zeta}v(s) d s\biggr\}\biggr]' F(\zeta) d\zeta = \phantom{aaaaaaaaaaaaaaaa}
$$
$$
\phantom{aaaaaaaaaaaaaaaa}=\exp\biggl\{-\il{t_0}{t}v(s) d s\biggr\} F(t) + \il{t_0}{t}\exp\biggl\{-\il{t_0}{\zeta}v(s) d s\biggr\}f(\zeta) d \zeta, \phh t \ge t_0.
$$
This together with (3.4) implies
$$
c + \il{t_0}{t}v(\zeta) d\zeta + F(t) \le \exp\biggl\{\il{t_0}{t}v(\zeta)d \zeta\biggr\}\biggl[c - \exp\biggl\{-\il{t_0}{t}v(s) d s\biggr\} F(t) + \biggr. \phantom{aaaaaaaaaaaaaaaaaaaaaaa}
$$
$$
\phantom{aaaaaaaaaaaaaaaaaaaaaaaaaaa} + \biggl.+\il{t_0}{t}\exp\biggl\{-\il{t_0}{\zeta}v(s) d s\biggr\}f(\zeta)d\zeta\biggr] + F(t), \phh t \ge t_0.
$$
Hence,
$$
c + \il{t_0}{t}v(\tau) d\tau + \il{t_0}{t}f(\tau) d\tau \le c\exp\biggl\{\il{t_0}{t}v(\zeta)d\zeta\biggr\} +
\il{t_0}{t}\exp\biggl\{\il{\zeta}{t}v(s) d s\biggr\} f(\zeta) d\zeta, \phh t \ge t_0.
$$
From here and from the condition (1.1) of the lemma it follows (1.2). The lemma is proved.

\vskip 10pt

{\bf 4.  Proof of Corollary 1.} If $t \ge \tau$ then by (1.3) we have
$$
u(t) \le u(\tau) + \il{\tau}{t}v(\zeta) u(\zeta) d\zeta + \il{\tau}{t} f(\zeta) d \zeta.
$$
Then in virtue of Lemma 2
$$
u(t) \le u(\tau)\exp\biggl\{\il{\tau}{t}v(s) d s\biggr\} + \il{\tau}{t}\exp\biggl\{\il{\zeta}{t}v(s) d s\biggr\} f(\zeta) d\zeta. \eqno (4.1)
$$
Let $t \le \tau$. Then by virtue of (1.3) we have
$$
u(t) \le u(\tau) + \il{t}{\tau}v(\zeta) u(\zeta) d\zeta + \il{t}{\tau} f(\zeta) d \zeta. \eqno (4.2)
$$
We set: $\widetilde{u}(t)\equiv u(-t), \ph \widetilde{v}(t)\equiv v(-t), \ph \widetilde{f}(t)\equiv f(-t), \ph  t \in (-b,-a)$. Rewrite (4.2) in the form
$$
\widetilde{u}(-t) \le \widetilde{u}(-\tau) + \il{t}{\tau}\widetilde{v}(-\zeta) \widetilde{u}(-\zeta) d\zeta + \il{t}{\tau} \widetilde{f}(-\zeta) d \zeta.
$$
Then
$$
\widetilde{u}(-t) \le \widetilde{u}(-\tau) + \il{-\tau}{-t}\widetilde{v}(\zeta) \widetilde{u}(\zeta) d\zeta + \il{-\tau}{-t} \widetilde{f}(\zeta) d \zeta, \ph a < t \le \tau < b.
$$
In this inequality we replace: $-t \to t, \ph -\tau \to \tau$. We obtain.
$$
\widetilde{u}(t) \le \widetilde{u}(\tau) + \il{\tau}{t}\widetilde{v}(\zeta) \widetilde{u}(\zeta) d\zeta + \il{\tau}{t} \widetilde{f}(\zeta) d \zeta, \ph -b < \tau \le t < -a.
$$
Then by already proven (4.1)  we will get
$$
\widetilde{u}(t) \le \widetilde{u}(\tau)\exp\biggl\{\il{\tau}{t}\widetilde{v}(s) d s\biggr\} + \il{\tau}{t}\exp\biggl\{\il{\zeta}{t}\widetilde{v}(s) d s\biggr\} \widetilde{f}(\zeta) d\zeta, \ph -b < \tau \le t < -a.
$$
It follows from here that
$$
\widetilde{u}(\tau)\exp\biggl\{\il{\tau}{t}\widetilde{v}(s) d s\biggr\} \ge \widetilde{u}(t) - \il{\tau}{t}\exp\biggl\{\il{\zeta}{t}\widetilde{v}(s) d s\biggr\} \widetilde{f}(\zeta) d\zeta, \ph -b < \tau \le t < -a.
$$
Multiplying both sides of this inequality by $\exp\biggl\{-\il{\tau}{t}\widetilde{v}(s)d s\biggr\}$ we obtain
$$
\widetilde{u}(\tau) \ge \widetilde{u}(t)\exp\biggl\{-\il{\tau}{t}\widetilde{v}(s)d s\biggr\} - \il{\tau}{t}\exp\biggl\{-\il{\tau}{\zeta}\widetilde{v}(s) d s\biggr\} \widetilde{f}(\zeta) d\zeta, \ph -b < \tau \le t < -a.
$$
It follows from here that
$$
u(-\tau) \ge u(-t) \exp\biggl\{\il{\tau}{t}v(-s) d(-s)\biggr\} + \il{\tau}{t}\exp\biggl\{\il{\tau}{\zeta}v(-s)d(-s)\biggr\}f(-\zeta)d(-\zeta),
$$
$-b < \tau \le t ,-a,$ or
$$
u(-\tau) \ge u(-t) \exp\biggl\{-\il{-t}{-\tau}v(s) d s\biggr\} - \il{-t}{-\tau}\exp\biggl\{\il{\zeta}{-t}v(s)d s\biggr\}f(\zeta)d \zeta, \ph -b < \tau \le t < -a.
$$
Replacing $-\tau$ by $t$ and $-t$ by $\tau$ from here we obtain
$$
u(t) \ge u(\tau) \exp\biggl\{-\il{\tau}{t}v(s) d s\biggr\} - \il{\tau}{t}\exp\biggl\{-\il{\zeta}{t}v(s) d s\biggr\} f(\zeta) d \zeta, \ph a < t \le \tau < b.
$$
This together with (4.1) implies (1.4). The corollary is proved.

\vskip 10pt

\centerline{\bf References}

\vskip 10pt

\noindent
1. L. Y. Adrianoba, Introduction to the theory of linear systems of differential equations.  \linebreak \phantom{a}
S. Peterburg, Publishers of St. Petersburg University, 1992.

\noindent
2. O. Akinyele, On Gronwall-Bellman-Bihary-Tipe Integral Inequalities in several \linebreak \phantom{a} variables with retardation. J. Math Anal. Appl., 104, 1--26 (1984).

\noindent
3. Bellman R. Stability theory of differential equations, Moscow, Izdatelstvo inostrannoj \linebreak \phantom{a}  literatury, 1954.

\noindent
4. Ch.-Ch. Yen and M.-H. Shih, The Gronwall-Bellman Inequality in Several variables. \linebreak \phantom{a} J. Math Anal. Appl., 86, 167--167 (1982).

\noindent
5. J. Chandra and B. A. Fleishman, On a Generalization of the Gronwall-Bellman Lemma \linebreak \phantom{a}  in Partially Ordered Banach Spaces. J. Math Anal. Appl., 31, 668--681 (1970).

\noindent
6. J. Conlan and Ch.-L. Wang, Gronwall-Bellman-type Inequalities for systems. J. Math.\linebreak \phantom{a}  Anal. Appl.,  120, 622--6030 (1986).

\noindent
7. Demidovich B. P.  Lectures on the mathematical  stability theory, Moscow, ''Nauka'', \linebreak \phantom{a}  1967.

\noindent
8. H. Fan, B. Zhang, Some New generalized Gronwall-Bellman Type Inequalities Arising \linebreak \phantom{a} In The Theory of Fractional Didderentia-integro Equations, WSEAS Transactions on \linebreak \phantom{a}  Mathematics, vol. 13, 2-14, pp. 820--829.

\noindent
9. G. A. Grigorian,  On two comparison tests for second-order linear  ordinary\linebreak \phantom{aa} differential equations (Russian) Differ. Uravn. 47 (2011), no. 9, 1225 - 1240; trans-\linebreak \phantom{aa} lation in Differ. Equ. 47 (2011), no. 9 1237 - 1252, 34C10.

\noindent
10. M. Hammy and M. A. Hammamy, Gronwall-Bellman type integral inequalities and \linebreak \phantom{a}  applications to global uniform stability. CUBO A Mathematical Journal, vol. 17, \linebreak \phantom{a}  N$^o$ 03 (53--70), 2015.

\noindent
11. L. Horvath, Gronwall-Bellman Type Integral Inequalities in Measure space. J. Math.\linebreak \phantom{a}  Anal. Appl., 202, 183--199 (1996).

\noindent
12. Z. A. Khan, Integral Inequalities of Gronwall-Bellman Type. Appl. Math., 2014, \linebreak \phantom{a} 5, 3484--3488.

\noindent
13. Z. A. Khan, Generalizations of Gronwall-Bellman Type Inequalities. Applied \linebreak \phantom{a} Mathematical Sciences, vol. 8, 2014, N$^o$ 57, 2817--1822.

\noindent
14. Z. Li and W.-SH. Wang, Some new nonlinear powered Gronwall-Bellman type retarded \linebreak \phantom{a} integral inequalities and their applications. J. Math. Ineq. vol. 13, Number 2 \linebreak \phantom{a} (2019), 553--564.

\noindent
15. L. Losonzi, A Generalization of The Gronwall-Bellman Lemma and its Applications. \linebreak \phantom{a} J. Math. Anal. Appl., 44, 701--709 (1973)

\noindent
16. Y. Louattassi, A new Generalization of Lemma Gronwall-Bellman. Applied Math. \linebreak \phantom{a} Sciences, vol. 6. 2012, N$^o$ 13, 621--628.

\noindent
17. M. Mekki, Boukerrioua, B. Kilani and M. L. Sahary, New explicit bounds on Gronwall- \linebreak \phantom{a} Bellman-Bihary-Gamedov integral inequalities and their weakly singular analogues with \linebreak \phantom{a} applications. Kraguevac J. Math., vol 44 (4) (2020) pp. 603--615.

\noindent
18. B. G. Pachpate, A Note on Gronwall-Bellman Identity. J. Math Anal. Appl., \linebreak \phantom{a} 44, 758--762 (1973).

\noindent
19.  B. G. Pachpate, On some generalizations of Bellman's Lemma. J. Math Anal. Appl.,\linebreak \phantom{a}  51, 141--150 (1975).

\noindent
20. S. Raffeq and S. Hussain, Anew Gronwall-Bellman type integral inequality and its \linebreak \phantom{a} application to fractional stochastic differential equations. Dynam. Syst. Appl., 28 \linebreak \phantom{a} N$^o$ 2 (2019), 259--273.

\noindent
21. James G. Taylor, Improved error bounds for the Liouville-Green (or WKB) \linebreak \phantom{a} Approximation. J. Math. Anal. Appl., 83, 79--89 (1982).

\noindent
22. G. Zhihong and L. Yongqing, Integral Inequalities of Gronwall-Bellman Type for \linebreak \phantom{a} Multy-distributions, J. Math Anal. Appl,, 183, 63--75 (1994).

\end{document}